\documentclass[oneside, a4paper,11pt,reqno]{amsart}
\usepackage{color}

\usepackage{amssymb,amsmath,amsthm}
\usepackage[T1]{fontenc}

\newtheorem{hypo}{Hypothesis}

\newtheorem{thm}[hypo]{Theorem}

\newtheorem{lem}[hypo]{Lemma}

\newtheorem{rqe}[hypo]{Remark}

\def\A{\mathcal{A}}

\def\PP{\mathbb{P}}
\def\RR{\mathbb{R}}
\def\ZZ{\mathbb{Z}}

\def\EE{\mathbb{E}}

\newcommand {\sur}[2] { \stackrel {\scriptstyle{#1}}{#2}}

\begin{document}
%%\maketitle 

 \vglue50pt

\centerline{ \Large \bf  The quenched limiting distributions }
\centerline{ \Large \bf  of  a  one-dimensional random walk in random scenery}

\bigskip

\medskip
 \renewcommand{\thefootnote}{\fnsymbol{footnote}}

 \centerline{Nadine Guillotin-Plantard\footnote{Institut Camille Jordan, Universit\'e   Lyon 1,  43, boulevard du 11 novembre 1918
69622 Villeurbanne.   Research partially  supported by ANR  (MEMEMOII)  2010 BLAN 0125 Email:  nadine.guillotin@univ-lyon1.fr}, 
Yueyun Hu\footnote{D\'epartement de Math\'ematiques (LAGA CNRS-UMR 7539) Universit\'e Paris 13, 99 avenue J.B. Cl\'ement, 93430  Villetaneuse.   Research partially  supported by ANR  (MEMEMOII)  2010 BLAN 0125 Email: yueyun@math.univ-paris13.fr}, and  Bruno Schapira\footnote{Centre de Math\'ematiques et Informatique, Aix-Marseille Universit\'e, 39, rue F. Joliot Curie, 13453 Marseille Cedex 13. Research partially  supported by ANR  (MEMEMOII)  2010 BLAN 0125. Email:  bruno.schapira@latp.univ-mrs.fr}}

\medskip

 \centerline{\it Universit\'e Lyon 1, Universit\'e Paris 13,  and Universit\'e Aix-Marseille}

% \medskip \centerline{\sl    This version:    08.07.2013}
\bigskip

{\leftskip=2truecm
\rightskip=2truecm
\baselineskip=15pt
\small

\noindent{\slshape\bfseries Summary.}        For a one-dimensional random walk in random scenery (RWRS) on $\ZZ$, we determine  its  quenched   weak limits   by applying  Strassen \cite{Strassen}'s functional law of the iterated logarithm.  As a consequence, conditioned on the random scenery, the one-dimensional RWRS does not converge in law, in contrast with the multi-dimensional case.

\medskip

 %\medskip

%%\centerline{This version: \today}
%\noindent{\slshape\bfseries 2000 Mathematics Subject Classification.} 60J80.

} %%%%%% End of narrower

\medskip

\section{Introduction}
Random walks in random sceneries were introduced independently by Kesten and Spitzer \cite{KS79} and by Borodin \cite{Bor79-1,Bor79-2}.
Let $S = (S_n)_{n \ge 0}$ be a random walk in $\mathbb{Z}^d$ starting at $0$,
i.e., $S_0 = 0$ and
$
\left(S_n-S_{n-1}\right)_{n \ge 1} \text{ is a sequence of i.i.d.\ } \mathbb{Z}^d \text{-valued random variables}.
$
Let $\xi = (\xi_x)_{x \in \mathbb{Z}^d}$ be a field of i.i.d.\ real random variables independent of $S$.
The field $\xi$ is called the random scenery.
The random walk in random scenery (RWRS) $K := (K_n)_{n \ge 0}$ is defined 
by setting $K_0 := 0$ and, for $n \in \mathbb{N}^{*}$,
\begin{equation}
K_n := \sum_{i=1}^n \xi_{S_i}.
\end{equation}
We will denote by $\mathbb{P}$ the joint law of $S$ and $\xi$. The law $\mathbb{P}$ is called the \emph{annealed} law, while the conditional law $\mathbb{P}(\cdot | \xi)$ is
called the \emph{quenched} law.

Limit theorems for RWRS have a long history, we refer to \cite{GuPo} or \cite{GuPoRe} for a complete review.  Distributional limit theorems for {\it quenched} sceneries (i.e. under the quenched law) are however quite recent. The first result in this direction that we are aware of was obtained by Ben Arous and \v{C}ern\'y \cite{BAC07}, in the case of a heavy-tailed scenery and planar random walk. In \cite{GuPo}, quenched central limit theorems (with the usual $\sqrt{n}$-scaling and Gaussian law in the limit) were proved 
for a large class of transient random walks.  More recently, in \cite{GuPoRe}, the case of the planar random walk was studied, the authors proved a quenched version of the annealed central limit theorem obtained by Bolthausen in \cite{Bo89}.

In this note we consider the case of the simple symmetric  random walk $(S_n)_{n\ge0}$ on $\ZZ$, the random scenery $(\xi_x)_{x\in \ZZ}$  is assumed to be centered with finite variance equal to one and there exists some $\delta>0$ such that $\EE(|\xi_0|^{2+\delta})<\infty$. 
We prove that under these assumptions, there is no quenched distributional limit theorem for $K$. In the sequel, for $-\infty \le a< b \le \infty$, we will denote by $\mathcal{AC}([a,b]\to {\mathbb R})$ the set of absolutely continuous functions defined on the interval $[a,b]$ with values in $\mathbb{R}$. Recall that  if $f\in\mathcal{AC}([a,b]\to {\mathbb R}) $, then the derivative of $f$ (denoted by $\dot{f} $) exists  almost everywhere and is  Lebesgue integrable  on
 $[a,b]$. Define  \begin{equation} \mathcal{K}^*  := \Big\{ f\in \mathcal{AC}( \RR \to \RR) :  f(0)=0, \int_{-\infty}^{\infty} (\dot{f}(x))^2 dx  \le 1\Big\}. \label{K*}  \end{equation}

\begin{thm}\label{T:0} 
For $\mathbb{P}$-a.e.\ $\xi$, under the quenched probability $\mathbb{P} \left(. \mid \xi \right)$, the process 
 $$\tilde{K}_n:= \frac{K_n} {(2 n^{3/2}Ê\log \log n)^{1/2}},  \qquad n >e^e, $$ 
does not converge in law. 
More precisely, for $\mathbb{P}$-a.e.\ $\xi$, under the quenched probability $\mathbb{P} \left(. \mid \xi \right)$, 
the limit points of the  law of 
$\tilde{K}_n,$ as $n \to \infty,  $
 under the topology of  weak convergence of measures,  are  equal to the set of the laws of random variables  in $\Theta_B $,  with  
 \begin{equation} \Theta_B  :=  \Big\{ \int _{-\infty}^{\infty} f(x) d L_1(x):    f \in  {\mathcal K}^*  \Big\},  \label{KB} \end{equation} where     $(L_1(x), x \in \RR)$ denotes the family of local times at time $1$ of a one-dimensional Brownian motion $B$ starting from $0$.
 
 The set $\Theta_B$ is closed for the topology of  weak convergence of measures, and is a compact subset of $L^2 ( (B_t)_{t\in[0,1]})$.
 
  \end{thm}

Let us mention that  the   set $\mathcal{K}^*$ directly comes from Strassen \cite{Strassen}'s limiting set.   The precise meaning of $ \int _{-\infty}^{\infty} f(x) d L_1(x)$ can be given by the integration by parts and the occupation times formula:    
\begin{equation}\label{occupation}
 \int _{-\infty}^{\infty}  f(x) d L_1(x)=- \int _{-\infty}^{\infty}  L_1(x) \dot {f}(x) dx = - \int_0^1 \dot{f}(B_s) ds, 
\end{equation}

\noindent where as before,  $\dot{f}$ denotes the almost everywhere derivative of $f$.

Instead of Theorem \ref{T:0}, we shall prove that there is no quenched limit theorem for the continuous analogue of $K$ introduced by Kesten and Spitzer \cite{KS79} and deduce  Theorem \ref{T:0}  by using a strong approximation for the one-dimensional RWRS. Let us define this continuous analogue: Assume that  $B:=(B(t))_{t\geq 0}$, $W:=(W(t))_{t\geq 0}$, $\tilde{W}:=(\tilde{W}(t))_{t\geq 0}$ are  three real Brownian motions starting from $0$,  defined on the same probability space and independent of each other.  For brevity, we shall write $W(x):= W(x)$ if $x \ge 0$ and $\tilde{W}(-x)$ if $x <0$ and say that $W$ is a two-sided Brownian motion.   We denote by $\mathbb{P}_{B}$, $\mathbb{P}_{W}$ the law of these processes.  We will also denote by $(L_t(x))_{t\geq 0,x\in\RR}$ a continuous version with compact support of the local time of the process $B$. 
We define the continuous version of the RWRS, also called {\it Brownian motion in Brownian scenery},  as 
$$Z_t:=\int_{0}^{+\infty}  L_{t} (x) d W(x) + \int_{0}^{+\infty} L_t(-x) d \tilde{W}(x) \equiv \int_{-\infty}^{+\infty} L_t(x) d W(x).$$
In dimension one, under the annealed measure,  Kesten and Spitzer \cite{KS79} proved that the process 
$(n^{-3/4} K([nt]) )_{t\geq 0}$  weakly converges in the space of continuous functions to the continuous process $Z= (Z(t))_{t\geq 0}$. Zhang \cite{Zh} (see also \cite{CKS99, KL98}) gave a stronger version of this result in the special case when the scenery has a finite moment  of order $2+\delta$ for some $\delta >0$, more precisely, there is a coupling of $\xi$, $S$, $B$
 and $W$ such that $(\xi, W)$ is independent of $(S,B)$ and for any $\varepsilon >0$, almost surely,
\begin{equation}\label{Zh}
\max_{0\leq m \leq n} | K(m) - Z(m)| = o(n ^{\frac{1}{2}+\frac{1}{2(2+\delta)} + \varepsilon} ) , \  \  n\rightarrow +\infty.
\end{equation}
Theorem \ref{T:0} will follow from this strong approximation and the following result.

\begin{thm}\label{T:1} $\PP_W$-almost surely,  under the quenched probability $\PP(\cdot \vert W)$,  the limit points of the  law of 
$$\tilde{Z}_t:= \frac{Z_t} {(2 t^{3/2}Ê\log \log t )^{1/2}},  \qquad t \to \infty,  $$
 under the topology of  weak convergence of measures,  are  equal to the set of the laws of random variables  in $\Theta_B $ defined in Theorem \ref{T:0}. Consequently under $\PP(\cdot \vert W)$, as $t \to \infty$,   $\tilde{Z}_t$ does not converge in law. 
 
 %, where 
 %\begin{equation} \label{KB} \mathcal{K}_B:=\Big\{ \int _{m_1}^{M_1} f(x) d L_1(x):  f\in \mathcal{C}([m_1,M_1] \to {\mathbb R} ) \mbox{ such that } f(0)=0,  \int_{m_1}^{M_1}  (\dot{f}(x))^2 dx \le 1\Big\}, \end{equation}
 %$m_1:= \inf_{0\le s \le 1} B_s$ and $M_1:= \sup_{0\le s \le 1} B_s$. 
 
 \end{thm}

To prove Theorem \ref{T:1}, we shall   apply  Strassen \cite{Strassen}'s functional law of the iterated logarithm applied to the two-sided  Brownian motion $W$; we shall also need to estimate the stochastic integral $\int g   (x) d L_1(x)$ for a   Borel function $g$,   see Section \ref{S:Proofs} for  the details.

\section{Proofs}\label{S:Proofs}

For a  two-sided one-dimensional Brownian motion  $(W(t), t \in {\mathbb R})$ starting from $0$, let us define for any $ \lambda > e^e$, 
$$ W_{\lambda} (t):= \frac{W(\lambda t)}{ (2 \lambda \log \log \lambda)^{1/2}},  \qquad   t \in \RR. $$

\begin{lem} \label{LIL} 

(i)   Almost surely, for any $s<0<r$ rational numbers, $(W_\lambda(t), s\le t \le r)$ is  relatively compact in the uniform topology and the set of its limit points is  $\mathcal{K}_{s,r}$, with  $$\mathcal{K}_{s,r}:=\Big\{ f\in \mathcal{AC}([s,r] \to {\mathbb R} ) :  f(0)=0, \int_s^{r} (\dot{f}(x))^2 dx  \le 1\Big\}.$$

(ii) There exists some finite random variable $\A_W$ only depending on $(W(x), x \in \RR)$ such that for all $\lambda \ge e^{36}$,  $$    \sup_{ t \in \RR, t  \neq0} { \vert W_\lambda(t)\vert \over \sqrt{ \vert t\vert  \log \log (  \vert t\vert + {1\over \vert t\vert } + 36)}} \le \A_W  < \infty .$$
\end{lem}

\begin{rqe} \label{R4} The statement  (i) is a reformulation of Strassen's theorem and holds in fact for all real numbers $s$ and $r$. Moreover, using the notation  $\mathcal{K}^*$ in \eqref{K*}, we  remark  that  $\mathcal{K}_{s,r}$ coincides with the restriction of $\mathcal{K}^*$ on  $[s, r]$:   for any $s<0<r$,  $$ {\mathcal K}_{s,r}= \Big\{ f_{\big \vert [s, r]}: f \in \mathcal{K}^*\Big\}.$$   \end{rqe}

{\noindent\bf Proof:}  (i) For any fixed $s<0<r$, by applying   Strassen's theorem (\cite{Strassen}) to the two-dimensional rescaled Brownian motion: $({W(\lambda r u) \over \sqrt{2\lambda r \log \log \lambda}}, {W(\lambda   s  u) \over \sqrt{2\lambda  \vert s\vert \log \log \lambda}})_{0\le u \le 1}$, we get that a.s.,  $(W_\lambda(t), s\le t \le r)$ is  relatively compact in the uniform topology  with $\mathcal{K}_{s, r}$ as  the set of   limit points.    By inverting a.s. and $s, r$, we obtain (i).

(ii)  By the classical law of the iterated logarithm for the Brownian motion $W$ (both at $0$ and at $\infty$), we get that $$ \widetilde \A_W := \sup_{x \in \RR, x \neq0} { \vert W(x) \vert \over  \sqrt{ \vert x\vert  \log \log (  \vert x\vert + {1\over \vert x\vert } + 36)}} $$ is a finite variable.  Observe  that  for any $t >0$ and $\lambda > e^{36}$, $ \log \log (  \lambda t  + {1\over  \lambda t } + 36)  \le \log \log \lambdaÊ+  \log \log (    t  + {1\over    t } + 36).$ The Lemma follows if we take for e.g.  $\A_W:= 2 \widetilde \A_W$. $\Box$

Next, we recall  some properties  of Brownian local times: The process $x \mapsto L_1(x)$ is a (continuous) semimartingale (by Perkins \cite{Perkins}), moreover, the quadratic variation of $x \mapsto L_1(x)$ equals $ 4 \int_{-\infty}^x L_1(z) d z$.  By Revuz and Yor (\cite{RevuzYor},   Exercise  VI (1.28)), for any locally bounded Borel function $f$, \begin{equation}\label{int-parts}  {1\over2} \int _{-\infty}^{\infty}  f(x) d L_1(x)= -  \int_0^{B_1} f(u) d u + \int_0^1 f(B_u) d B_u.\end{equation}

Let us define for all $\lambda > e^e$ and $ n\ge 0$, $$ H_\lambda:= \int _{-\infty}^{\infty} W_\lambda(x) d L_1(x), \qquad H_\lambda^{(n)}:= \int_{-n}^n W_\lambda(x) d L_1(x),$$

\noindent with $H^{(0)}_\lambda=0$.  Denote by $\EE_B$ the expectation with respect to the law of $B$.

\begin{lem}\label{L9}  There exists some positive constant $c_1$ such that for any $\lambda> e^{36}$ and $n \ge0$, we have  \begin{eqnarray}    \EE_B \Big\vert   H_\lambda- H_\lambda^{(n)} \Big\vert &\le&  c_1\,  e^{- {n^2\over4}}\, \A_W, \label{h2nb} 
  \\ \EE_B  \Big( \int _{-\infty}^{\infty}  f(x) d L_1(x) \Big)^2   & \le&  16\,  s(f), \label{fL(x)}   
    \\ \EE_B \Big \vert \int _{-\infty}^{\infty}  f(x) d L_1(x)- \int _{-n}^{n}  f(x) d L_1(x)   \Big\vert & \le&   4 \,  \sqrt{2 s(f)} \, e^{-{n^2\over4}}, \label{fL(x)n} 
\end{eqnarray}   for any   Borel function $f: \RR \to \RR$ such that $ s(f):=  \sup_{0\le u \le 1} \EE_B \Big[ f^2(B_u) \Big]  < \infty.$ 
\end{lem}

Remark that if $f$ is bounded, then $s(f) \le \sup_{x \in \RR} f^2(x)$.

{\noindent\bf Proof:}  We first prove that there exists some positive constant $c_2$ such that for all $n \ge0$ and $\lambda> e^{36}$,   \begin{equation}\label{hl2n}  \EE_B \Big[ (H_\lambda- H_\lambda^{(n)})^2\Big] \le  c_2 \, \A_W^2. \end{equation}

\noindent
 In fact, by applying  \eqref{int-parts} and  using  the Brownian isometry  for $f(x)= W_\lambda(x) 1_{( \vert x \vert >n)} $,  we get that $$ \EE_B \Big[ (H_\lambda- H_\lambda^{(n)})^2\Big] \le 8  \EE_B \Big[ F_{n, \lambda}(B_1)^2\Big]  + 8  \EE_B \Big[  \int_0^1 (W_\lambda(B_u))^2 1_{(\vert B_u \vert >n)} d u \Big],$$ 

\noindent with $F_{n,\lambda}(x):= \int_0^x W_\lambda(y) 1_{(\vert y \vert >n)} dy $ for any $x \in \RR$. By Lemma \ref{LIL} (ii), $$ \vert F_{n,\lambda}(x) \vert \le  \A_W \, \Big\vert \int_0^x  (\vert y \vert \log \log ( \vert y \vert + {1\over \vert y \vert } +36))^{1/2} d y \Big \vert \le c_3\, \A_W \, (1+x^2), \qquad \forall x \in \RR,$$ with some constant $c_3>0$.  Hence $\EE_B \Big[ F_{n, \lambda}(B_1)^2\Big]   \le 6\, c_3^2\, \A_W^2$.    In the same way, $\EE_B\big[ (W_\lambda(B_u))^2\big]  \le \A_W^2\, \EE_B \big[ \vert B_u \vert \log \log ( \vert B_u\vert + {1\over \vert B_u \vert}+36)\big]$ which is   integrable for $u \in (0, 1]$. Then \eqref{hl2n} follows.

To check \eqref{h2nb}, we remark that $ H_\lambda- H_\lambda^{(n)}=0$ if $\sup_{0\le u \le 1} \vert B_u \vert \le n$. Then by Cauchy-Schwarz' inequality and \eqref{hl2n}, we have that \begin{eqnarray*}
\EE_B \Big\vert   H_\lambda- H_\lambda^{(n)} \Big\vert &=& \EE_B \Big[ \big\vert   H_\lambda- H_\lambda^{(n)} \big\vert \, 1_{( \sup_{0\le u \le 1} \vert B_u \vert >n)} \Big] \\
	 & \le& \sqrt{\EE_B \Big[ (H_\lambda- H_\lambda^{(n)})^2\Big] } \, \sqrt{ \PP_B\Big( \sup_{0\le u \le 1} \vert B_u \vert >n\Big)} \\
	 &\le& \sqrt{c_2} \, \A_W\, \sqrt{2}\, e^{- {n^2\over4}},
\end{eqnarray*}

\noindent by the standard Gaussian tail: $\PP_B\big( \sup_{0\le u \le 1} \vert B_u \vert > x\big) \le 2 e^{- x^2/2}$ for any $x>0$. Then    we get \eqref{h2nb}. 

To prove \eqref{fL(x)}, we use again   \eqref{int-parts} and the Brownian isometry to arrive at  $$  \EE_B   \Big( \int _{-\infty}^{\infty}  f(x) d L_1(x) \Big)^2   \le 8 \EE_B \Big[ G^2(B_1) \Big]+   8 \int_0^1 \EE_B \Big[ f^2(B_u)\Big] d u \le 8 \EE_B \Big[ G^2(B_1) \Big]+ 8 s(f),$$

\noindent with $G(x):= \int_0^x f(y) dy$ for any $x \in \RR$. By Cauchy-Schwarz' inequality, $(G(x))^2 \le \Big\vert x  \int_0^x   f^2(y)  dy \Big\vert$ for any $x \in\RR$, from which we use  the integration by parts for the density of $B_1$ and deduce that $\EE_B \Big[ G^2(B_1) \Big] \le \EE_B \Big[ f^2(B_1) \Big] $.  Then  \eqref{fL(x)} follows.

  Finally for  \eqref{fL(x)n}, we use   \eqref{fL(x)} to see that $$\EE_B \Big ( \int _{-\infty}^{\infty}  f(x) d L_1(x)- \int _{-n}^{n}  f(x) d L_1(x)   \Big)^2 =  \EE_B \Big ( \int _{-\infty}^{\infty}  f(x) 1_{(\vert x\vert >n)} d L_1(x)   \Big)^2 \le 16 s(f), $$ for any $n$.  Then  \eqref{fL(x)n} follows from the Cauchy-Schwarz inequality and the Gaussian tail, exactly in the same way as   \eqref{h2nb}.   $\Box$

Recalling    \eqref{KB} for the definition of $\Theta_B$. For any $p>0$,  it is easy to see that $\Theta_B \subset L^p(B)$, since from Cauchy-Schwarz' inequality, using the relation (\ref{occupation}), we deduce that 
$$\big( \int _{-\infty}^{\infty}  f(x) d L_1(x)\big)^2 \le \big(  \int_{-\infty}^{\infty} (L_1(x))^2 dx \big) \big(  \int_{-\infty}^{\infty}  (\dot{f}(x))^2 dx\big) \le \sup_{x} L_1(x) \in L^p(B),$$

\noindent see Cs\'aki \cite{Csaki89}, Lemma 1 for the tail of $\sup_{x} L_1(x) $.  Write $d_{L^1(B)}(\xi, \eta)$ for the distance in $L^1(B)$ for any $\xi, \eta \in L^1(B)$.

\begin{lem}  \label{L:conv} $\PP_W$-almost surely, $$d_{L^1(B)}( H_\lambda, \Theta_B) \to 0, \qquad \mbox{ as } \lambda \to \infty,$$ where  $\Theta_B$ is defined in \eqref{KB}.    Moreover, $\PP_W$-almost surely for any $\xi \in  \Theta_B$,  $ \liminf_{\lambda \to \infty} d_{L^1(B)}( H_\lambda, \xi)=0.$
\end{lem}

%%As before, for any $ f\in \mathcal{C}([m_1,M_1] \to {\mathbb R} )$  such that $ f(0)=0$  and  $ \int_{m_1}^{M_1}  (\dot{f}(x))^2 dx \le 1$,   the integral $\int _{m_1}^{M_1} f(x) d L_1(x)$ should be understood as  $\int _{m_1}^{M_1}  f(x) d L_1(x) =  - \int_0^1 {\dot f} (B_s) ds  .$  

{\noindent\bf Proof:}  Let $\varepsilon>0$.  Choose a large $n=n(\varepsilon)$ such that $c_1 e^{-n^2/4} \le \varepsilon$.  By Lemma \ref{LIL} (i),  for all large $\lambda \ge \lambda_0(W, \varepsilon, n)$, there exists some function $g =g_{\lambda, W, \varepsilon, n}\in {\mathcal K}_{-n, n}$ such that $\sup_{ \vert x \vert \le n} \vert W_\lambda(x) - g(x) \vert \le \varepsilon$.   Applying  \eqref{fL(x)}  to $ f(x)= (W_\lambda(x)- g(x)) 1_{(\vert x\vert \le n)}$ which is a bounded   by $\varepsilon$,  we get that $$ \EE_B \Big\vert H_\lambda^{(n)} - \int_{-n}^n g(x) d L_1(x) \Big\vert \le 4\sqrt{s(f) }\le 4 \varepsilon.$$

We extend $g$ to $\RR$ by letting $g(x)= g(n)$ if $x\ge n$ and $g(x)= g(-n)$ if $x \le -n$, then $g \in \mathcal K^*$ and $ \int_{-\infty}^\infty g(x) d L_1(x)=  \int_{-n}^n g(x) d L_1(x)$.   By the triangular inequality and  \eqref{h2nb},  
$$ \EE_B \Big\vert H_\lambda  - \int_{-\infty}^\infty g(x) d L_1(x) \Big\vert \le 4 \varepsilon+  \EE_B \Big\vert H_\lambda- H^{(n)}_\lambda\Big\vert \le (4 + c_1 \A_W) \varepsilon. $$

It follows that $d_{L^1(B)}( H_\lambda, \Theta_B) \le (4 + c_1 \A_W) \varepsilon$. Hence $\PP_W$-a.s., $\limsup_{\lambda\to\infty} d_{L^1(B)}( H_\lambda, \Theta_B) \le (4 +c_1 \A_W) \varepsilon$, showing the first part in the Lemma.

For the other part of the Lemma, let $h \in \mathcal K^*$ such that $\xi=  \int _{-\infty}^{\infty}  h(x) d L_1(x)$.  Observe that $\vert h(x) \vert \le \sqrt{\Big\vert x  \int_0^x (\dot{h}(y))^2 dy \Big\vert} \le \sqrt{\vert x\vert}$ for all $x \in \RR$, $s(h)= \sup_{0\le u \le 1} \EE_B[ h^2(B_u)] \le \EE_B[ \vert B_1\vert]$, then for any $\varepsilon>0$,  we may use \eqref{fL(x)n} and choose an integer $n=n(\varepsilon)$ such that $(c_1 + 4\sqrt{2}) e^{-n^2/4} \le \varepsilon$ and $$ d_{L^1(B)}( \xi, \xi_n)   \le \varepsilon,$$

\noindent where $\xi_n:= \int _{-n}^{n}  h(x) d L_1(x)$.  Applying Lemma \ref{LIL} (i) to the restriction of $h$ on $[-n, n]$, we may find  a sequence $\lambda_j= \lambda_j(\varepsilon, W, n) \to \infty$ such that $\sup_{ \vert x \vert \le n} \vert W_{\lambda_j}(x) - h(x) \vert \le \varepsilon$.  By applying \eqref{fL(x)} to $f(x)= (W_{\lambda_j}(x) - h(x)) 1_{(\vert x\vert \le n)}$, we have that $$ d_{L^1(B)}( H^{(n)}_{\lambda_j}, \, \xi_n) \le 4 \varepsilon.$$

By \eqref{h2nb} and the choice of $n$, $ d_{L^1(B)}( H^{(n)}_{\lambda_j}, H_{\lambda_j}) \le \varepsilon \A_W$ for all large $\lambda_j$, it follows from the triangular inequality that $$ d_{L^1(B)} ( \xi, H_{\lambda_j}) \le (5 + \A_W )\varepsilon,$$
implying that $\PP_W$-a.s., $ \liminf_{\lambda \to \infty} d_{L^1(B)}( H_\lambda, \xi) \le (5+ \A_W )\varepsilon \to 0$ as $\varepsilon\to0$.   $\Box$
  
 We now are ready to give the proof of Theorems \ref{T:1} and \ref{T:0}.
 
 {\noindent\bf Proof of Theorem \ref{T:1}.} Firstly, we remark that by Brownian scaling, $\PP_W$-a.s.,  \begin{equation}\label{Z}
\frac{Z_t}{t^{3/4}} \sur{(d)}=  -  \int_{m_1}^{M_1} \frac{1}{t^{1/4}} W( \sqrt{t} y) d L_1(y) .  \end{equation} 

 In fact,  by the change of variables $x= y\sqrt{t}$, we get 
 $$\int_{-\infty}^{+\infty}  L_{t} (x) d W(x) = \sqrt{t} \int_{-\infty}^{+\infty}  \left( \frac{L_{t}( y \sqrt{t})}{\sqrt{t}}\right) dW( y \sqrt{t})$$
which has the same distribution as
$$\sqrt{t} \int_{-\infty}^{+\infty} L_{1}( y)  dW( y \sqrt{t})$$
from  the scaling property of the local time of the Brownian motion. 
Since $(L_1(x))_{x \in \RR}$ is a continuous semi-martingale, independent from the process $W$, from the formula of integration by parts, we get that $\mathbb{P}_W$ -a.s.,
$$\sqrt{t} \int_{-\infty}^{+\infty} L_{1}( y)  dW( y \sqrt{t})= - t^{3/4} \int_{m_1}^{M_1} \left( \frac{W(\sqrt{t} y)} {t^{1/4}} \right)  d L_1(y),$$

\noindent yielding \eqref{Z}.  The first part of Theorem \ref{T:1} follows from Lemma \ref{L:conv}. 

Let $(\zeta_n)_n$ be a sequence of random variables in $\Theta_B$, each $\zeta_n$ being associated to a function $f_n\in \mathcal{K}^*$. The sequence of the (almost everywhere) derivatives of $f_n$ is then a bounded sequence in the Hilbert space $L^2(\RR)$, so we can extract a subsequence which weakly converges. Using the definition of the weak convergence and the relation $(\ref{occupation})$, $(\zeta_n)_n$ converges almost surely and the closure of $\Theta_B$ follows. Since the sequence $(\zeta_n)_n$ is bounded in $L^p(B)$ for any $p\geq 1$, the convergence also holds in $L^2(B)$. Therefore $\Theta_B$ is a compact set of $L^2(B)$ as closed and bounded subset. $\Box$

{\noindent\bf Proof of Theorem \ref{T:0}.}
We use the strong approximation of Zhang \cite{Zh} : there exists on a suitably enlarged probability space,  a coupling of $\xi$, $S$, $B$
and $W$ such that $(\xi, W)$ is independent of $(S,B)$ and for any $\varepsilon >0$, almost surely,
$$
\max_{0\leq m \leq n} | K(m) - Z(m)| = o(n ^{\frac{1}{2}+\frac{1}{2(2+\delta)} + \varepsilon} ) , \  \  n\rightarrow +\infty.
$$
From the independence of $(\xi, W)$ and $(S,B)$, we deduce that for $\PP$-a.e. $(\xi, W)$, under the quenched probability $\PP(.\vert \xi,W)$, the limit points of the 
 laws of $\tilde{K}_n$ and $\tilde{Z}_n$ are the same ones. 
Now, by adapting the proof of Theorem \ref{T:1}, we have that for $\PP$-a.e. $(\xi, W)$, under the quenched probability $\PP(.\vert \xi,W)$, the limit points of the 
 laws of  $\tilde{Z}_n$, as $n\rightarrow \infty$,  under the topology of  weak convergence of measures,  are  equal to the set of the laws of random variables  in 
 $\Theta_B $. It gives that for $\PP$-a.e. $(\xi, W)$, under the quenched probability $\PP(.\vert \xi,W)$, the limit points of the 
 laws of  $\tilde{K}_n$, as $n\rightarrow \infty$,  under the topology of  weak convergence of measures,  are  equal to the set of the laws of random variables  in 
 $\Theta_B $ and Theorem \ref{T:0} follows.
$\Box$

       \medskip
    
{\noindent\bf Acknowledgments.} We are grateful to Mikhail Lifshits  for  interesting discussions. 
The authors thank the referee for recommending various improvements in exposition.


\begin{thebibliography}{13}
\baselineskip=10pt

\bibitem{BAC07}
Ben Arous, G. and \v{C}ern\'y, J., (2007) Scaling limit for trap models on $\mathbb{Z}^d$,
{\it Ann. Probab.} 35 (6), 2356 -- 2384.


\bibitem{Bo89} Bolthausen, E. (1989)
A central limit theorem for two-dimensional random walks in random sceneries,
{\it Ann.\ Probab.} 17, 108--115.

\bibitem{Bor79-1} Borodin, A.N. (1979)
A limit theorem for sums of independent random variables defined on a recurrent random walk,
{\it Dokl. Akad. nauk SSSR} 246 (4), 786 -- 787.

\bibitem{Bor79-2} Borodin, A.N. (1979) 
Limit theorems for sums of independent random variables defined in a transient random walk,
in Investigations in the Theory of Probability Distributions, IV, Zap, Nauchn. Sem. Leningrad. Otdel. Mat. Inst. Steklov. (LOMI) 85, 17-29. 237 244.

\bibitem{Csaki89} Cs\'aki, E. (1989) An integral test for the supremum
of Wiener local time.  {\it Probab. Th. Rel. Fields} 83, 207--217.

\bibitem{CKS99} Cs\'aki, E.,  K\"onig, W. and Shi, Z. (1999)
An embedding for the Kesten-Spitzer random walk in random scenery,
{\it Stochastic Process. Appl.} 82 (2), 283-292.

\bibitem{GuPo}
Guillotin-Plantard, N. and  Poisat, J. (2013)
Quenched central limit theorems for random walks in random scenery,
{\it Stochastic Process. Appl.} 123 (4), 1348 -- 1367. 

\bibitem{GuPoRe}
Guillotin-Plantard, N., Poisat, J. and Dos Santos, R.S. (2013)
A quenched central limit theorem for planar random walks in random sceneries,
{\it Submitted}. 


% \bibitem{JeulinYor}   Jeulin, T. and  Yor, M.  (1977) {\it Grossissement d'une filtration et semi-martingales: formules explicites.} S\'eminaire de Probabilit\'es, XII,  78 -- 97, Lecture Notes in Math., 649, Springer, Berlin, 1978.

\bibitem{KS79}
Kesten, H. and Spitzer, F. (1979)
A limit theorem related to a new class of self-similar processes,
{\it Z. Wahrsch. Verw. Gebiete} 50 (1), 5--25.


\bibitem{KL98} Khoshnevisan, D. and Lewis, T.M. (1998)
A law of the iterated logarithm for stable processes in random scenery,
{\it Stochastic Process. Appl.} 74 (1), 89--121.


%\bibitem{MansuyYor}  Mansuy, R. and  Yor, M. (2006)  {\it Random times and enlargements of filtrations in a Brownian setting. } Lecture Notes in Mathematics, 1873. Springer-Verlag, Berlin. 


\bibitem{Perkins}  Perkins, E.  (1982).  Local time is a semimartingale. {\it Z. Wahrsch. Verw. Gebiete } {\bf 60}   no. 1, 79 -- 117. 

\bibitem{RevuzYor}  Revuz, D. and Yor, M. (1999). {\it Continuous martingales and Brownian motion.} Third edition.  Springer-Verlag, Berlin. 


\bibitem{Strassen} Strassen, V. (1964). An invariance principle for the law of the iterated logarithm.  {\it Z. Wahrsch. Verw. Gebiete }  {\bf 3},  Issue 3, 211 -- 226. 

\bibitem{Zh} Zhang, L. (2001) The strong approximation for the Kesten-Spitzer random walk. {\it Statistics \& Probability Letters} {\bf 53}, 21 -- 26.

\end{thebibliography}
\end{document}